\documentstyle[amsmath,amssymb,12pt]{article}

\newcommand{\rom}[1]{{\rm #1}}

\allowdisplaybreaks

\begin{document}

\setcounter{page}{1} \setcounter{section}{0} \thispagestyle{empty}

\newtheorem{definition}{Definition}
\newtheorem{remark}{Remark}
\newtheorem{proposition}{Proposition}
\newtheorem{theorem}{Theorem}
\newtheorem{corollary}{Corollary}
\newtheorem{lemma}{Lemma}

\newcommand{\indlim}{\operatornamewithlimits{ind\,lim}}
\newcommand{\Ffin}{{\cal F}_{\mathrm fin}}
\newcommand{\Fext}{{\cal F}_{\mathrm ext}}
\newcommand{\D}{{\cal D}}
\newcommand{\N}{{\Bbb N}}
\newcommand{\C}{{\Bbb C}}
\newcommand{\Z}{{\Bbb Z}}
\newcommand{\R}{{\Bbb R}}
\newcommand{\Rp}{{\R_+}}
\newcommand{\eps}{\varepsilon}
\newcommand{\supp}{\operatorname{supp}}
\newcommand{\la}{\langle}
\newcommand{\ra}{\rangle}
\newcommand{\const}{\operatorname{const}}
\renewcommand{\emptyset}{\varnothing}
\newcommand{\di}{\partial}
\newcommand{\hotimes}{\hat\otimes}

\renewcommand{\author}[1]{\medskip{\Large #1}\par\medskip}

\newcommand{\pii}{\pi_{\nu\otimes\sigma}}
\newcommand{\RR}{{\cal R}}
\newcommand{\RX}{{\RR\times X}}
\newcommand{\ZZ}{\Z_{+,\,0}^\infty}

\begin{center}{\Large \bf  The square of white noise as a Jacobi field }\end{center}

\author{Eugene Lytvynov}

\noindent{\sl Department of Mathematics\\ University of Wales Swansea\\ Singleton Park\\ Swansea
SA2 8PP\\
United Kingdom 
\\[2mm]
 {\rm E-mail: e.lytvynov@swansea.ac.uk}}

\begin{abstract}

\noindent We identify the representation of the square of white
noise obtained by L.~Accardi, U.~Franz and M.~Skeide in [{\it Comm.\ Math.\ Phys.}\ {\bf 228} (2002),
123--150] with the Jacobi field of a L\'evy process of Meixner's
type.

\end{abstract}


\section{Formulation of the result}

The problem of developing a stochastic calculus for higher powers
of white noise, i.e., ``nonlinear stochastic  calculus'', was
first stated by Accardi, Lu, and Volovich in \cite{1}. Since the
white noise is an operator-valued distribution, in order to solve
this problem one needs an appropriate renormalization procedure.
In \cite{2,3}, it was proposed to renormalize the commutation
relations and then to look for Hilbert space  representations of
them. Let us shortly discuss this approach.

We will use $\R^d$, $d\in\N$,  as  an underlying
space.
Let ${b}(x)$, $x\in \R^d$, be an operator-valued distribution
satisfying the canonical commutation relations:
\begin{gather} [{b}(x),{b}(y)] =
[{b}^\dag(x),{b}^\dag(y)]=\pmb0 ,\notag\\ [{b}(x),{b}^\dag(y)]
=\delta(x-y)\pmb1.\label{z7eawr76}\end{gather}
Here, $[A,B]{:=}AB-BA$ and ${b}^\dag(x)$ is the dual operator of
${b}(x)$. Denote \begin{equation} B_x{:=}{b}(x)^2,\quad
B_x^\dag{:=}{b}^\dag(x)^2,\quad N_x{:=}{b}^\dag(x){b}(x),\quad x\in
\R^d.\label{sdz}\end{equation} One wishes to derive from \eqref{z7eawr76} the commutation
relations satisfied by the operators $B_x,B_x^\dag, N_x$.
To this end, one needs to make sense  of the square of the
delta function, $\delta(x)^2$. But it is known from the
distribution theory that \begin{equation}\delta(x)^2= c
\delta(x),\label{bhefbhf}\end{equation}
where $c\in\C$ is an arbitrary constant (see \cite{2} for a justification of this formula
and bibliographical references).

Thus, using \eqref{z7eawr76} and   formula \eqref{bhefbhf} as a renormalization, we get
\begin{gather} [B_x,B_y^\dag]=2c\delta(x-y)
\pmb1+4\delta(x-y)N_y,\notag\\ [N_x,B_y^\dag]=2\delta(x-y)B_y^\dag,\notag\\
[N_x,B_y]=-2\delta(x-y)B_y,\notag\\
[N_x,N_y]=[B_x,B_y]=[B_x^\dag,B_y^\dag]=\pmb0\label{huuh}\end{gather}
(see\ \cite[Lemma~2.1]{4}).

Let ${\cal S}(\R^d)$ denote the Schwartz space of rapidly
decreasing functions on $\R^d$.
 For each $\varphi\in {\cal S}(\R^d)$, we
introduce \begin{equation}\label{dfhuio}
B(\varphi){:=}\int_{\R^d} \varphi(x) B_x\, dx,\quad B^\dag(\varphi){:=}\int_{\R^d} \varphi(x) B^ \dag_x\,
dx,\quad
N(\varphi){:=}\int_{\R^d}\varphi(x)N_x\,dx. \end{equation} By \eqref{huuh},
\begin{gather} [B(\varphi),B^\dag(\psi)]=2c\la \varphi,\psi\ra
\pmb1+4N(\varphi\psi),\notag\\ [N(\varphi),B^\dag(\psi)]=2B^\dag(\varphi\psi),\notag\\
[N(\varphi),B(\psi)]=-2B(\varphi\psi),\notag\\
[N(\varphi),N(\psi)]=[B(\varphi),B(\psi)]=[B^\dag(\varphi),B^\dag(\psi)]=\pmb0,\qquad \phi,\psi
\in {\cal S}(\R^d).\label{jwaui}\end{gather}
Here, $\la\cdot,\cdot\ra$ denotes the scalar product in
$L^2(\R^d,dx)$.
The Lie algebra with generators $B(\varphi),
B^\dag(\varphi),N(\varphi)$, $\varphi\in {\cal S}(\R^d)$,
and a central element $\pmb 1$ with relations \eqref{jwaui} is called
the square of white noise (SWN) algebra.

Now, one is interested in a
Hilbert space representation of the SWN algebra with a cyclic vector
$\Phi$ satisfying $B(\varphi)\Phi=0$ (which is called a Fock
representation). In \cite{2}, it was shown that a Fock
representation of the SWN algebra exists if and only if the
constant $c$ is strictly positive. In what follows,  we will suppose, for simplicity of notations
 that  $c=2$.

Let us now recall the Fock representation of the SWN algebra
constructed in \cite{5} (see also references therein).

For a real separable Hilbert space $\cal H$, denote by ${\cal F}({\cal
H})$ the symmetric Fock space over $\cal H$:
$$ {\cal F}({\cal H})=\bigoplus_{n=0}^\infty {\cal H}^{\hotimes
n}n!,$$ where $\hotimes $ stands for the symmetric tensor product.
Thus, each $f\in{\cal F}({\cal H})$ is of the form
$f=(f^{(n)})_{n=0}^\infty$, where $f^{(n)}\in{\cal H}^{\hotimes n}$ and $\|f\|_{{\cal F}({\cal H})}^2=\sum_{n=0}
^\infty \|f^{(n)}\|_{{\cal H}^{\hotimes n}}^2n!$\,.
Now take $\cal H$ to be $L^2(\R^d,dx)\otimes \ell_2$, where the
$\ell_2$ space
has the orthonormal basis $(e_n)_{n=1}^\infty$, $e_n=(0,\dots,0,\underbrace{1}_{\text {$n$th
place}},0,\dots)$.

Denote by $\frak F$ the linear subspace of ${\cal F}(L^2(\R^d,dx)\otimes \ell_2)$
that is the linear span of the
vacuum vector $\Omega=(1,0,0,\dots)$ and vectors of the form $(\varphi\otimes \xi)^{\otimes
n}$, where  $\varphi\in{\cal S}(\R^d)$ and  $\xi\in \ell_{2,0}$, $n\in\N$.
Here, $\ell_{2,0}$ denotes the linear subspace of $\ell_2$ consisting
of finite  vectors, i.e., vectors  of the form
$\xi=(\xi_1,\xi_2,\dots,\xi_m,0,0,\dots)$, $m\in\N$.  The set  $\frak F$ is evidently a dense
subset of
${\cal F}(L^2(\R^d,dx)\otimes \ell_2)$.

Denote by $J^+,J^0,J^-$ the linear operators in $\ell_2$ with domain
$\ell_{2,0}$ defined by the following formulas:
\begin{align}J^+e_n&=\sqrt{n(n+1)}\,e_{n+1},\notag\\
J^0e_n&=n e_n,\notag\\
J^-e_n&=\sqrt{(n-1)n}\,e_{n-1},\qquad n\in\N.\label{vefug}\end{align}

Now, for each $\varphi,\psi\in{\cal S}(\R^d)$ and $\xi\in\ell_{2,0}$, we set
\begin{align} B^\dag (\varphi)(\psi\otimes\xi)^{\otimes n}&=2(\varphi\otimes
e_1)\hotimes (\psi\otimes\xi)^{\otimes n}+2n((\varphi\psi)\otimes(J^+\xi))^{\otimes
n},\notag\\ N(\varphi)(\psi\otimes \xi)^{\otimes n}&= 2n((\varphi\psi)\otimes
J^0\xi)^{\otimes n},\notag\\  B(\varphi)(\psi\otimes\xi)^{\otimes n}&= 2n\la
\varphi,\psi\ra \xi_1(\psi\otimes \xi)^{\otimes(n-1)}+2n ((\varphi\psi)\otimes(J^-\xi))^{\otimes
n},\label{hgfaduz}
\end{align}
where $n\in\N$, and $(\psi\otimes \xi)^ {\otimes 0}{:=}\Omega$. Thus, \begin{align} B^\dag (\varphi)&=2 A^+(\varphi\otimes
e_1)+2A^0(\varphi\otimes J^+),\notag\\ N(\varphi)&=2A^0(\varphi\otimes
J^0),\notag\\ B(\varphi)&=2A^-(\varphi\otimes e_1)+2A^0(\varphi\otimes
J^-),\label{suivdu}\end{align} where $A^+(\cdot)$, $A^0(\cdot)$, and $ A^-(\cdot)$ are the creation, neutral, and
annihilation operators in ${\cal F}(L^2(\R^d,dx)\otimes \ell_2)$,
respectively. The operator $B^\dag(\varphi)$ is the restriction of
the adjoint operator of $B(\varphi)$ to $\frak F$, while the operator
$N(\varphi)$ is Hermitian. It is easy to see that the operators $B^\dag(\varphi), N(\varphi),
B(\varphi)$, $\varphi\in{\cal S}(\R^d)$, constitute a representation of the SWN
algebra.

In what follows, the closure of a closable operator $A$  will be
denoted by $\widetilde A$. Since the adjoint operators of $B^\dag(\varphi)$, $N(\varphi)$,
$B(\varphi)$ are densely defined, they are closable.

The last part of \cite{5} is devoted to studying 
those classical infinitely divisible processes which are built
from the SWN in a similar way as the Wiener and Poisson processes
are built from the usual white noise. So, for each parameter $\beta\ge
0$, we define \begin{equation}\label{usdi} X_\beta(x){:=}B_x^\dag+B_x+\beta N_x,\qquad
x\in\R^d.\end{equation}  Notice that we want a formally self-adjoint process,
so the parameter $\beta$ must be real (we also exclude from consideration the case $\beta<0$,
since it may be treated by a trivial transformation of the case
$\beta>0$).

In view of \eqref{z7eawr76} and \eqref{sdz}, the only privileged parameter is
$\beta=2$, when $X_\beta(x)$  becomes the renormalized square of
the classical white noise ${b}^\dag(x)+{b}(x)$, see \cite[Section~3]{4}.

Analogously to \eqref{dfhuio}, we introduce, for each $\varphi\in{\cal
S}(\R^d)$, \begin{equation} X_\beta(\varphi){:=}\int_{\R^d}\varphi(x)X_\beta(x)\,
dx=B^\dag(\varphi)+B(\varphi)+\beta N(\varphi).\label{zsguow}\end{equation}
As easily seen,  $\widetilde X_\beta(\varphi )$ is a
self-adjoint operator.

 In the case $d=1$, it was shown in \cite{5} that the quantum
 process $(\widetilde X_\beta(\chi_{[0,t]}))_{t\ge0}$ ($\chi_\Delta$ denoting the indicator function of a set $\Delta$)
 is associated with a classical  L\'evy process
 $(Y_\beta(t))_{t\ge0}$, which is a gamma process for $\beta=2$,
 a Pascal process for $\beta>2$, and a Meixner process for
 $0\le\beta<2$. (One has, of course, to extend the SWN algebra in order to include
 the operators indexed by the indicator functions, for example, to take the set $L^2(\R,dx)\cap L^\infty (\R,dx)$
 instead of ${\cal S}(\R)$.)

We also refer to \cite{4,AB,5} and references therein for a discussion of
other aspects of the SWN.

On the other hand, in papers \cite{KL,Ly3,Ly,belyme} (see also
\cite{silva,beme1,bere3,beme3}), the Jacobi field of the L\'evy processes of  Meixner's type, i.e., the gamma,
Pascal, and Meixner processes, was studied. Let us shortly explain
this approach.

Let ${\cal S}'(\R^d)$ be the Schwartz space of tempered distributions. The ${\cal S}'(\R^d)$ is the
dual space of ${\cal S}(\R^d)$ and the
dualization between ${\cal S}'(\R^d)$ and ${\cal S}(\R^d)$ is
given by the scalar product in $L^2(\R^d,dx)$. We will
preserve
the symbol $\la\cdot,\cdot\ra$ for this dualization. Let ${\cal C}({\cal S}'(\R^d))$
denote the cylinder $\sigma$-algebra on ${\cal S}'(\R^d)$.

For
each $\beta\ge 0$, we define a probability measure $\mu_\beta$ on
$( {\cal S}'(\R^d)),{\cal C}({\cal S}'(\R^d))$ by its Fourier
transform
\begin{equation}\label{rew4w}
\int_{{\cal S}'(\R^d)} e^{i\la \omega,\varphi\ra}\,
\mu_\beta(d\omega)=\exp\bigg[\int_{\R\times \R^d}
(e^{is\varphi(x)}-1-is\varphi(x))\,\nu_\beta(ds)\,dx\bigg],\qquad
\varphi\in {\cal S}(\R^d),\end{equation}
where the measure $\nu_\beta$ on $\R$ is specified as
follows.

Let $\tilde\nu_\beta$ denote the
probability measure on $(\R,{\cal B}(\R))$ whose orthogonal
polynomials $(\widetilde P_{\beta,\, n})_{n=0}^\infty$ with
leading coefficient 1 satisfy the recurrence relation
\begin{gather}\label{ghfd}s\widetilde P_{\beta,\, n}(s)=
\widetilde P_{\beta,\, n+1}(s)+\beta(n+1)\widetilde
P_{\beta,\, n}(s)+n(n+1)\widetilde P_{\beta,\, n-1}(s),\\
n\in\Z_+,\, \widetilde P_{\beta,\,-1}(s){:=}0.\notag\end{gather}
By  \cite[Ch.~VI, sect.~3]{Chihara},
$(\widetilde P_{\beta,\, n})_{n=0}^\infty$ is a system of
polynomials of Meixner's type, the measure $\widetilde\nu_\beta$
is uniquely determined by the above condition and is given as
follows. For $\beta\in[0,2)$, $$ \tilde\nu_\beta(ds)=
\frac{\sqrt{4-\beta^2}}{2\pi}\, \big|\Gamma\big(1+i
(4-\beta^2)^{-1/2}s\big)\big|^2\,\exp\big[-s2(4-\beta^2)^{-1/2}\arctan
\big(\beta(4-\beta^2)^{-1/2}\big) \big]\,ds $$
($\tilde\nu_\beta$ is a Meixner distribution), for $\beta=2$
$$\tilde\nu_2(ds)=\chi_{(0,\infty)} (s)e^{-s}s\,ds$$
($\tilde\nu_2$ is a gamma distribution), and for $\beta>2$
$$\tilde\nu_\beta(ds)=(\beta^2-4)\sum_{k=1}^\infty
p_\beta^k\,k\,\delta_{\sqrt{\beta^2-4}\,k},\qquad
p_\beta{:=}\frac{\beta-\sqrt{\beta^2-4}}{\beta+\sqrt{\beta^2-4}}$$
($\tilde\nu_\beta$ is now a Pascal distribution).

Notice that, for each $\beta\ge0$, $\tilde\nu(\{0\})=0$, and
hence, we may define \begin{equation}
\nu_\beta(ds){:=}\frac1{s^2}\,\tilde\nu_\beta(ds).\label{iwegf}\end{equation}
Then, $\mu_\beta$  is the measure of gamma
noise for $\beta=2$, Pascal noise for $\beta>2$, and Meixner noise
for $\beta\in[0,2)$. Indeed, for each $\beta\ge0$, $\mu_\beta$ is
a generalized process on ${\cal S}'(\R^d)$ with independent values
(cf.~\cite{GV}). Next, for each $\varphi\in{\cal S}(\R^d)$, we
have \begin{equation}\int_{{\cal S}'(\R^d)}\la
\omega,\varphi\ra^2\,\mu_\beta(d\omega)=\int_{\R^d}\varphi(x)^2\,
dx.\label{uieriaev}\end{equation} Hence, for each $f\in L^2(\R^d,dx)$, we may
define, in a standard way, the random variable $\la\cdot,f\ra$
from $L^2({\cal S}'(\R^d),d\mu_\beta)$ satisfying \eqref{uieriaev}
with $\varphi=f$.

Then, for each open, bounded  set $\Delta\subset \R^d$, the
distribution $\mu_{\beta,\Delta}$ of the random variable
$\la\cdot,\chi_\Delta\ra$ under $\mu_\beta$ is given as
follows. For $\beta>2$\rom, $\mu_{\beta,\Delta}$ is
the negative binomial \rom(Pascal\rom) distribution $$
\mu_{\beta,\Delta}=(1-p_{\beta})^{|\Delta|}\sum_{k=0}^\infty
 \frac{\big(|\Delta|\big)_k}{k!}\, p_\beta^k \, \delta_{\sqrt{\beta^2-4}\,k-2|\Delta|/(\beta+\sqrt{\beta^2-4})}, $$ where
  for $r> 0$ $(r)_0{:=}1$, $(r)_k{:=}r(r+1)\dotsm (r+k-1)$\rom,
$k\in\N$\rom. For $\beta=2$, $\mu_{2,\Delta}$ is the Gamma
distribution
$$\mu_{2,\Delta}(ds)=\frac{(s+|\Delta|)^{|\Delta|-1}e^{-(s+|\Delta|)}}{\Gamma(|\Delta|)}\,\chi_{(0,\infty)}(s+|\Delta|)\,
ds.$$ Finally\rom, for $\beta\in[0,2)$\rom,
\begin{multline*}\mu_{\beta,\Delta}(ds)=\frac{(4-\beta^2)^{(|\Delta|-1)/2}}{2\pi
\Gamma(|\Delta|) }\big|
\Gamma\big(|\Delta|/2+i(4-\beta^2)^{-1/2}(s+\beta|\Delta|/2\big)
\big|^2\\ \times\exp\big[
-(2s+\beta|\Delta|)(4-\beta^2)^{-1/2}\arctan\big(\beta(4-\beta^2)^{-1/2}\big)
\big]\,ds.\end{multline*}
Here, $|\Delta|{:=}\int_\Delta dx$.

We denote by ${\cal P}({\cal S}'(\R^d))$ the set of continuous
polynomials on ${\cal S}'(\R^d)$, i.e., functions on ${\cal S}'(\R^d)$ of the
form $$F(\omega)=\sum_{i=0}^n\la\omega^{\otimes
i},f^{(i)}\ra,\qquad\omega^{\otimes 0}{:=}1,\ f^{(i)}\in{\cal
S}(\R^d)^{\hotimes i},\ i=0,\dots,n,\ n\in\Z_+. $$ The greatest number
$i$ for which $f^{(i)}\ne0$ is called the power of a polynomial.
 We denote by ${\cal P}_n({\cal S}'(\R^d))$ the set of
continuous polynomials of power $\le n$.

The set  ${\cal P}({\cal S}'(\R^d))$ is a dense subset
of $ L^2({\cal S}'(\R^d), d\mu_\beta)$. Let ${\cal P}^\sim
_n({\cal S}'(\R^d))$ denote the closure of ${\cal P}_n({\cal S}'(\R^d))$ in
$L^2({\cal S}'(\R^d),d\mu_\beta )$,  let ${\bf P}_n({\cal S}'(\R^d))$,
$n\in\N$, denote the orthogonal difference ${\cal P }^\sim_n({\cal
S}'(\R^d))\ominus{\cal P}^\sim_{n-1}({\cal S}'(\R^d))$, and let ${\bf
P}_0({\cal S}'(\R^d)){:=}{\cal P }^\sim_0({\cal S}'(\R^d))$.
We evidently have the orthogonal
decomposition \begin{equation}\label{tswes} L^2({\cal
S}'(\R^d),d\mu_\beta )=\bigoplus_{n=0}^\infty{\bf P}_n({\cal
S}'(\R^d)).
\end{equation}

For a monomial $\la\omega^{\otimes n},f^{(n)}\ra$, $f^{(n)}\in{\cal
S}(\R^d)^{\hotimes n}$, we denote by ${:}\la \omega^{\otimes
n},f^{(n)}\ra{:}$ the orthogonal projection of $\la\omega^{\otimes
n},f^{(n)}\ra$ onto ${\bf P}_n({\cal S}'(\R^d))$. The set $\{{:}\la \omega^{\otimes
n},f^{(n)}\ra{:},\ f^{(n)}\in{\cal S}(\R^d)^{\hotimes n}\}$ is dense in ${\bf
P}_n({\cal S}'(\R^d))$.

Denote by $\ZZ $ the set of all sequences $\alpha$ of the form
$\alpha=(\alpha_1,\alpha_2,\dots,\alpha_n,0,0,\dots)$,
$\alpha_i\in\Z_+$, $n\in\N$. Let $|\alpha|{:=}\sum_{i=1}^\infty
\alpha_i$, evidently $|\alpha|\in\Z_+$.
For each $\alpha\in\ZZ$, $1\alpha_1+2\alpha_2+\dots=n$, $n\in\N$,
and for any function $f^{(n)}:(\R^d)^n\to\R$ we define a function $D_\alpha
f^{(n)}:(\R^d)^{|\alpha|}\to\R$ by setting \begin{align*}(D_\alpha
f^{(n)})(x_1,\dots,x_{|\alpha|}){:=}& f^{(n)}(x_1,\dots,x_{\alpha_1},
\underbrace{x_{\alpha_1+1},x_{\alpha_1+1}}_{\text{2 times }},
\underbrace{x_{\alpha_1+2},x_{\alpha_1+2}}_{\text{2 times
}},\dots,
\underbrace{x_{\alpha_1+\alpha_2},x_{\alpha_1+\alpha_2}}_{\text{2
times }},\notag\\ &\quad
\underbrace{x_{\alpha_1+\alpha_2+1},x_{\alpha_1+\alpha_2+1},,x_{\alpha_1+\alpha_2+1}}_{\text{3
times }},\dots).\end{align*}
We define a scalar product on ${\cal S}(\R^d)^{\hotimes n}$
by setting for any $f^{(n)},g^{(n)}\in{\cal S}(\R^d)^{\hotimes n}$
\begin{gather}
(f^{(n)},g^{(n)})_{{\cal F}^{(n)}_{\mathrm Ext}(L^2(\R^d,dx))}
{:=}\sum_{\alpha\in\ZZ:\, 1\alpha_1+2\alpha_2+\dots=n}K_\alpha
\int_{X^{|\alpha|}}(D_\alpha f^{(n)})(x_1,\dots,x_{|\alpha|})\notag
\\ \times (D_\alpha g^{(n)})(x_1,\dots,x_{|\alpha|}) \,dx_1\dotsm dx_{|\alpha|},\label{1.5}\end{gather}
where \begin{equation}\label{t7rr5r5} K_\alpha=\frac{n!}{\alpha_1!\, 1^{\alpha_1}\alpha_2!\,
2^{\alpha_2}\dotsm}\,.
\end{equation}
Let
${\cal F}_{{\mathrm Ext}}^{(n)}(L^2(\R^d,dx))$ be
the closure of ${\cal S}(\R^d)^{\hotimes n}$ in the norm generated by \eqref{1.5}, \eqref{t7rr5r5}.
The extended Fock space $
{\cal F}_{{\mathrm Ext}}(L^2(\R^d,dx))$
over $L^2(\R^d,dx)$ is defined as
\begin{equation}\label{1.5a}
{\cal F}_{{\mathrm Ext}}(L^2(\R^d,dx)){:=}\bigoplus_{n=0}^\infty {\cal F}^{(n)}
_{{\mathrm Ext}}(L^2(\R^d,dx))\,n!,\end{equation}
where ${\cal F}^{(n)}
_{{\mathrm Ext}}(L^2(\R^d,dx)){:=}\R$. We also denote by $\Omega$
the vacuum vector in ${\cal F}_{{\mathrm Ext}}(L^2(\R^d,dx))$:
$\Omega=(1,0,0,\dots)$.

For any $f^{(n)}, g^{(n)}\in{\cal S}(\R^d)^{\hat\otimes n}$,
$n\in\N$, we have \begin{equation}\label{uigsd} \int_{{\cal S}'(\R^d)} {:}\la\omega^{\otimes
n},f^{(n)}\ra{:}\, {:}\la\omega^{\otimes
n},g^{(n)}\ra{:}\,\mu_\beta(d\omega)=  (f^{(n)},g^{(n)})_{{\cal F}^{(n)}
_{{\mathrm Ext}}(L^2(\R^d,dx))}\,n!\,. \end{equation}
Therefore, for each $f^{(n)}\in {\cal F}^{(n)}
_{{\mathrm Ext}}(L^2(\R^d,dx))$, we can define,  a random variable ${:}\la\cdot^{\otimes
n},f^{(n)}\ra{:}$ from $L^2({\cal S}'(\R^d),d\mu_\beta)$ such that
equality \eqref{uigsd} remains true for any \linebreak $f^{(n)},g^{(n)}\in {\cal F}^{(n)}
_{{\mathrm Ext}}(L^2(\R^d,dx))$, and furthermore  \begin{multline}\label{tsetse} {\cal F}_{{\mathrm
Ext}}(L^2(\R^d,dx))\ni f=(f^{(n)})_{n=0}^\infty\mapsto \\ \mapsto
U_\beta f=(U_\beta f)(\omega)=\sum_{n=0}^\infty {:}\la\omega^{\otimes
n},f^{(n)}\ra{:} \in L^2({\cal S}'(\R^d),d\mu_\beta)\end{multline} is unitary.

We denote by ${\cal F}_{\mathrm fin }({\cal S}(\R^d))$ the dense subset
of ${\cal F}_{{\mathrm
Ext}}(L^2(\R^d,dx))$ consisting of vectors of the form
$(f^{(0)},f^{(1)},\dots,f^{(n)},0,0,\dots)$, where
$f^{(i)}\in{\cal S}(\R^d)^{\hat\otimes i}$.  For each $\beta\ge0$
and each $\varphi\in{\cal S}(\R^d)$, we define an operator
$a_\beta(\varphi)$ on ${\cal F}_{\mathrm fin }({\cal S}(\R^d))$ by
the following formula:
$$a_\beta(\varphi)=a^+(\varphi)+\beta
a^0(\varphi)+a^-(\varphi).$$ Here\rom, $a^+(\xi)$ is the
standard creation operator\rom: \begin{equation}
a^+(\varphi)f^{(n)}{:=}\varphi\hotimes f_n,\qquad f^{(n)}\in {\cal S}(\R^d)^{\hotimes n},\ n\in\Z_+,
\label{udf}\end{equation} $a^0(\varphi)$ is the
standard neutral operator\rom: \begin{equation}
(a^0(\varphi)f^{(n)})(x_1,\dots,x_n)=\big(\varphi(x_1)+\dots+\varphi(x_n)\big)f_n(x_1,\dots,x_n),
\label{uwewe}\end{equation} and
\begin{equation}a^-(\varphi)=a^-_1(\varphi)+a^-_2(\varphi),\label{wuef}\end{equation} where
$a^-_1(\varphi)$ is the standard annihilation operator\rom:
\begin{equation}
(a^-_1(\varphi)f^{(n)})(x_1,\dots,x_{n-1})=n\int_{\R^d} \varphi(x)
f^{(n)}(x,x_1,\dots,x_{n-1})\,dx, \label{uzgwe}\end{equation} and
\begin{equation}
(a^-_2(\varphi)f^{(n)})(x_1,\dots,x_{n-1})=n(n-1)(\varphi(x_1)
f^{(n)}(x_1,x_1,x_2,x_3,\dots,x_{n-1}))^{\sim} ,
\label{gwvdgw}\end{equation}  $(\cdot)^\sim$ denoting symmetrization of a function.

Denote by $\di_x^\dag$, $\di_x$ the standard creation and
annihilation operators at point $x\in\R^d$: $$ \di^\dag_x
f^{(n)}=\delta_x\hotimes f^{(n)},\quad \di_x f^{(n)}(x_1,\dots,x_{n-1})=nf^{(n)}(x,x_1,\dots,x_{n-1}).$$
Then, at least formally, we have the following representation:
\begin{equation} a^+(\varphi)=\int_{\R^d}\varphi(x)\di_x^\dag\,
dx, \quad a^0(\varphi)=\int_{\R^d}\varphi(x)\di_x^\dag\di_x\,
dx,\quad a^-(\varphi)=\int_{\R^d}
\varphi(x)(\di_x+\di^\dag_x\di_x^2)\,dx,
\label{isuzg}\end{equation} so that \begin{equation}
a_\beta(\varphi)=\int_{\R^d}\varphi(x)
(\di^\dag_x+\beta\di^\dag_x\di_x+\di_x+\di^\dag_x\di_x^2)\,dx.
\label{ned456}\end{equation} (In fact, equalities  \eqref{isuzg}, \eqref{ned456} may be given
a precise meaning, cf.~\cite{KL,Ly3}.)

The operators $a_\beta(\varphi)$, $\varphi\in{\cal S}(\R^d)$,
 are essentially self-adjoint on ${\cal F}_{\mathrm fin }({\cal
 S}(\R^d))$ and the image of any $\tilde a_\beta(\varphi)$, $\varphi\in{\cal
S}(\R^d)$, under the unitary $U_\beta$ is the operator of
multiplication by the random variable $\la\cdot,\varphi\ra$. Thus,
$(\tilde a(\varphi))_{\varphi\in{\cal S}(\R^d)}$ is the Jacobi
field of  $\mu_\beta$, see   \cite{bere,berre,Ly8,belyme} and the references therein.

The functional realization of the operators $a^+(\varphi)$,
$a^0(\varphi)$, $a^-(\varphi)$, i.e., the explicit action of the the image of these operators under
the unitary $U_\beta $ is discussed in
\cite{KL,Ly3}.

A direct computation shows that the operators
$2a^+(\varphi),2a^0(\varphi),2a^-(\varphi)$,
$\varphi\in{\cal S}(\R^d)$, satisfy the commutation relations
\eqref{jwaui}, and hence generate a SWN algebra.
In fact, we have the following result:

\begin{theorem} For each $\beta\ge0$, there exists a unitary
operator $$I_\beta: {\cal F}(L^2(\R^d,dx)\otimes\ell_2)\to {\cal
F}_{\mathrm Ext} (\L^2(\R^d,dx))$$ such that
$I_\beta\Omega=\Omega$ and the operators $\widetilde
X_\beta(\varphi)$\rom, $\widetilde B^\dag(\varphi)$\rom,
$\widetilde N(\varphi)$\rom,  $\widetilde B(\varphi)$\rom, $\varphi\in{\cal
S}(\R^d)$\rom,  are
unitarily isomorphic under $I_\beta$ to two times the operators
$\tilde a(\varphi)$\rom, $\tilde a^+(\varphi)$\rom, $\tilde
a^0(\varphi)$\rom, $\tilde a^-(\varphi)$, respectivlely\rom.

\end{theorem}

Notice that  the unitary operator $$ {\cal
U}_\beta{:=}U_\beta I_\beta: {\cal F}(L^2(\R^d,dx)\otimes
\ell_2)\to L^2({\cal S}'(\R^d),d\mu_\beta) $$ has the following
properties\rom: ${\cal U}_\beta\Omega=1$ and $$ {\cal U}_\beta
\widetilde X_\beta(\varphi){\cal
U}_\beta^{-1}=2\la\cdot,\varphi\ra\cdot\, ,\qquad \varphi\in{\cal
S}(\R^d)$$
(compare with \cite{5})

By virtue of \eqref{dfhuio}, \eqref{usdi}, \eqref{isuzg}, and
\eqref{ned456}, we  get from Theorem~1:
\begin{equation} B_x=2(\di_x+\di_x^\dag\di_x^2),\quad 
N_x=2\di^\dag_x\di_x,\quad
B_x^\dag=2\di_x^\dag, \quad\label{qwegz9}\end{equation} and $$ X_\beta(x)=
2(\di_x^\dag+\beta\di_x^\dag\di_x+\di_x+\di^\dag_x\di_x^2),\qquad
x\in\R^d\label{z8we7we8}$$
(where the equalities are to be understood in the sense of the unitary
isomorphism). The reader is advised to compare \eqref{qwegz9} with
the informal representation \eqref{sdz}.

\section{Proof of the theorem}

The proof of Theorem~1 is essentially based on the results of
\cite{Ly}. By \eqref{suivdu} and \eqref{zsguow}, we get, for each
$\varphi\in{\cal S}(\R^d)$, $$ X_\beta(\varphi)=2(A^+(\varphi\otimes e_1)+A^0(\varphi\otimes J_\beta)
+A^-(\varphi\otimes e_1)),$$ where $$ J_\beta{:=}J^++\beta
J^0+J^-.$$ By \eqref{vefug}, the operator $J_\beta$ is given by a
Jacobi matrix (see e.g.\ \cite{b}). Furthermore, $J_\beta$ is
essentially self-adjoint on $\ell_{2,0}$
and, by \eqref{ghfd}, $\tilde \nu_\beta$ is the spectral measure
of $\widetilde J_\beta$. The latter means that there exists a
unitary operator $$ I_\beta^{(1)}:\ell_2\to
L^2(\R,d\tilde\nu_\beta)$$ such that $I_\beta^{(1)}e_1=1$ and, under $I_\beta^{(1)}$, the
operator $\widetilde J_\beta$ goes over into the operator of
multiplication by $s$.

Next, by \eqref{iwegf}, the operator $$ L^2(\R,d\tilde\nu_\beta)\ni
f\mapsto I_\beta^{(2)}f=(I_\beta^{(2)}f)(s){:=}f(s)s\in
L^2(\R,d\nu_\beta)$$ is unitary.
Setting $$I_\beta^{(3)}{:=}I_\beta^{(2)}I_\beta^{(1)}:\ell_2\to L^2(\R,d\nu_\beta),$$ we get a
unitary operator such that $I_\beta^{(3)}e_1=(I_\beta^{(3)}e_1)(s)=s$ and, under $I_\beta^{(3)}$, $\widetilde
J_\beta$ goes over into the operator of multiplication by $s$.

Using $I_\beta^{(3)}$, we can naturally construct a unitary
operator $$
I_\beta^{(4)}: {\cal F}(L^2(\R^d,dx)\otimes\ell_2)\to {\cal F}(L^2(\R^d,dx)\otimes
L^2(\R,d\nu_\beta))$$ such that $I_\beta^{(4)}\Omega=\Omega$ and,
under $I_\beta^{(4)}$, the operator $X_\beta(\varphi)$ goes over
into the operator $$ {\cal X}_\beta(\varphi)=2(A^+(\varphi\otimes s)+A^0(\varphi\otimes s)+A^-(\varphi
\otimes s)).$$

It follows from \cite{Ly} that there exists a unitary operator $$
I_\beta^{(5)}:{\cal F}(L^2(\R^d,dx)\otimes
L^2(\R,d\nu_\beta))\to L^2({\cal S}'(\R^d),d\mu_\beta)$$ such
that $I_\beta^{(5)}\Omega=1$ and, under $I_\beta^{(5)}$, the
operator $\widetilde{\cal X}_\beta(\varphi)$ goes over into the
operator of multiplication by $2\la\cdot,\varphi\ra$.

We define the unitary $$
I_\beta{:=}U_\beta^{-1}I_{\beta}^{(5)}I_\beta^{(4)}:{\cal F}(L^2(
\R^d,dx)\otimes\ell_2)\to {\cal F}_{\mathrm Ext}(L^2(\R^d,dx)),$$
where $U_\beta$ is given by \eqref{tsetse}. We evidently  get $I_\beta\Omega=\Omega$ and $\tilde a(\varphi)=I_{\beta}^{-1}\widetilde X_\beta(\varphi)
I_\beta^{-1}$, $\varphi\in{\cal S}(\R^d)$.

Next, we denote by $\frak G$ the subset of ${\cal F}_{\mathrm Ext}(L^2(\R^d,dx))$
defined as the linear span of $\Omega$ and the
vectors of the form $\varphi^{\otimes n}$, where $\varphi\in{\cal
S}(\R^d)$ and  $n\in\N$.  We note: $$ (I_\beta^{(3)}e_n)(s)
=P_{\beta,n}(s), \qquad n\in\N
,$$ where  $$P_{\beta,n}(s){:=}s\widetilde P_{\beta,n-1}(s),\qquad n\in\N, $$ and 
$(\widetilde P_{\beta,n})_{n=0}^\infty$ are defined by
\eqref{ghfd}. Hence, by \cite[Sect.~4 and Corollary~5.1]{Ly}, $$ {\frak
G}\subset I_\beta{\frak F}.$$ Furthermore, by \eqref{vefug},
\eqref{hgfaduz}, \eqref{udf}--\eqref{gwvdgw} and by
\cite[Corollary~5.1]{Ly}, we get: \begin{align} I_\beta B^\dag(\varphi)
I_\beta^{-1}\restriction {\frak G}&= a^+(\varphi)\restriction {\frak G},\notag\\
I_\beta N(\varphi)
I_\beta^{-1}\restriction {\frak G}&= a^0(\varphi)\restriction {\frak G},\notag\\
I_\beta B(\varphi)
I_\beta^{-1}\restriction {\frak G}&= a^-(\varphi)\restriction {\frak G},\qquad \varphi\in{\cal
S}(\R^d).\label{uiqweui}\end{align}

We now endow ${\cal F}_{\mathrm fin}({\cal S}(\R^d))$ with the
topology of the topological direct sum of the spaces ${\cal F}_n({\cal
S}(\R^d))$. Thus, the convergence in ${\cal F}_{\mathrm fin}({\cal
S}(\R^d))$ means the uniform finiteness and the coordinate-wise
convergence in each ${\cal F}_n({\cal
S}(\R^d))$. As easily seen, $\frak G$ is a dense subset of ${\cal F}_{\mathrm fin}({\cal
S}(\R^d))$. Since the operators $a^+(\varphi)$, $a^0(\varphi)$, and
$a^-(\varphi)$ act continuously on ${\cal F}_{\mathrm fin}({\cal S}(\R^d))$
and since ${\cal F}_{\mathrm fin}({\cal S}(\R^d))$ is continuously
embedded into ${\cal F}_{\mathrm Ext}(L^2(\R^d,dx))$
(cf.~\cite[p.~37]{KL}), the closure of the operators $a^+(\varphi)$, $a^0(\varphi)$, and
$a^-(\varphi)$  restricted to $\frak G$ coincides with $\tilde a^+(\varphi)$, $\tilde a^0(\varphi)$, and
$\tilde a^-(\varphi)$, respectively. Hence, by \eqref{uiqweui}, $\widetilde
B^\dag(\varphi)$, $\widetilde N(\varphi)$, and $\widetilde
N(\varphi)$ are extensions of the operators $\tilde a^+(\varphi)$, $\tilde a^0(\varphi)$, and
$\tilde a^-(\varphi)$, respectively.

Finally, analogously to the proof of \cite[Theorem~6.1]{Ly}, we
conclude that $I_\beta{\frak F}$ is a subset of the domain of $\tilde a^+(\varphi)$, respectively
$\tilde a^0(\varphi)$, respectively $\tilde a^-(\varphi)$, and
furthermore
\begin{align*} I_\beta B^\dag (\varphi) I_\beta^{-1}&=\tilde
a^+(\varphi)\restriction I_\beta{\frak F},\\
I_\beta N (\varphi) I_\beta^{-1}&=\tilde
a^0(\varphi)\restriction I_\beta{\frak F},\\
I_\beta B(\varphi) I_\beta^{-1}&=\tilde
a^-(\varphi)\restriction I_\beta{\frak F},\qquad \varphi\in{\cal
S}(\R^d).\end{align*}
This yields:
\begin{align*} I_\beta \widetilde B^\dag (\varphi) I_\beta^{-1}&=\tilde
a^+(\varphi),\\
I_\beta \widetilde N (\varphi) I_\beta^{-1}&=\tilde
a^0(\varphi),\\
I_\beta\widetilde  B(\varphi) I_\beta^{-1}&=\tilde
a^-(\varphi),\qquad \varphi\in{\cal
S}(\R^d),\end{align*}
which concludes the proof.
\vspace{2mm}

\noindent {\bf Acknowledgements.} I am grateful to L.~Accardi for
his inspiring me to consider this problem. I would also like to
thank U. Franz for interesting discussions.

\end{document}